\documentclass[12pt,a4paper]{article}
\usepackage{inputenc, amssymb, amsmath, amsthm}
\usepackage[english]{babel}
\makeatletter
\thm@style{plain}

\def\th@plain{%
  \thm@headfont{\bfseries}%
  \itshape 
  \thm@notefont{\rm}%
}
\def\thm@indent{\hspace*{\parindent}}
\makeatother

\newcommand{\flr}[1]{\left\lfloor #1 \right\rfloor}

\def\({\left(}
\def\){\right)}
\newcommand{\al}{\alpha}
\newcommand{\bt}{\beta}
\newcommand{\de}{\delta}
\newcommand{\bde}{\overline{\delta}}
\newcommand{\e}{\varepsilon}

\newcommand{\BA}{\boldsymbol{A}}

\newcommand{\baa}{\overline{\vphantom{\beta}\alpha}} 
\newcommand{\bb}{\overline{\beta}} 

\newcommand{\be}{\begin{equation}}
\newcommand{\ee}{\end{equation}}

\let\epsilon\varepsilon
\let\phi\varphi
\let\le\leqslant

\def\gcd{\mathop{\ensuremath{\text{\textup{\tt gcd}}}}}
\newtheorem{theorem}{Theorem}
\newcommand{\dokvo}{{\it Proof.} }
\newcommand\inte{\int\limits}

\begin{document}

\centerline{\bf RECOVERING FOURIER COEFFICIENTS}
\centerline{\bf OF MODULAR FORMS}
\centerline{\bf AND FACTORING OF INTEGERS}

\medskip

\centerline{Sergei~N.~Preobrazhenski\u i\footnote[1]{%
{\it Preobrazhenskii Sergei Nikolayevich} ---
Department of Mathematical Analysis, Faculty of Mechanics and Mathematics,
Lomonosov Moscow State University.}}

\bigskip

\hbox to \textwidth{\hfil\parbox{0.9\textwidth}{%
\small It is shown that if a function
defined on the segment $[-1,1]$ has sufficiently
good ap\-pro\-xi\-ma\-ti\-on by partial sums of the Legendre polynomial expansion,
then, given the function's Fourier coefficients $c_n$ for some subset of $n\in[n_1,n_2]$,
one can approximately recover them for all $n\in[n_1,n_2]$.
As an application, a new approach to factoring of integers is given.

\medskip

}\hfil}

\bigskip


{\bf 1. Introduction.}
Let a function $f(x)$
being a polynomial of degree $K$ in the interval $\left[-\frac12,\frac12\right)$
be continued periodically to the entire real axis with period~$1$:
\[
f(x)=b_0+b_1x+\ldots+b_Kx^K,\quad x\in\left[-\frac12,\frac12\right),\quad
b_K\ne0,\quad f(x+1)=f(x).
\]
A question is this: Can we find all Fourier coefficients of the
function $f(x)$ if we know $K+1$ values of the coefficients
$c_{p_0}$, $\ldots$, $c_{p_K}$? Obviously, one may set up
the system of linear equations
\begin{equation*}
\begin{pmatrix}
a_{1,1}, & \dots, & a_{1,K+1} \\
a_{2,1}, & \dots, & a_{2,K+1} \\
\hdotsfor{3} \\
a_{K+1,1}, & \dots, & a_{K+1,K+1}
\end{pmatrix}
\begin{pmatrix}
b_0 \\
b_1 \\
\vdots \\
b_K
\end{pmatrix}
=
\begin{pmatrix}
c_{p_0} \\
c_{p_1} \\
\vdots \\
c_{p_K}
\end{pmatrix}
.
\end{equation*}
So if we prove that the matrix of the system is nonsingular
then we can uniquely identify the polynomial and find the remaining Fourier coefficients
of the function $f(x)$.

In this note we show that if a function
defined on the segment $[-1,1]$ has sufficiently
good ap\-pro\-xi\-ma\-ti\-on by partial sums of the Legendre polynomial expansion,
then, given the function's Fourier coefficients $c_n$ for some subset of $n\in[n_1,n_2]$,
one can approximately recover them for all $n\in[n_1,n_2]$.
As an application, a new approach to factoring of integers is given.
Among factoring algorithms now in use the fastest ones are:
the elliptic curve method, the quadratic sieve and the number field sieve method.
The quadratic sieve algorithm was introduced by Pomerance~[1, 2] in 1981.
The heuristic complexity of the algorithm is
$$
\exp\left((1+o(1))(\log n)^{1/2}(\log\log n)^{1/2}\right)
$$
arithmetic operations, where $n$ is the number to be factored.
The elliptic curve method of H.~Lenstra~[3, 4] appeared in 1986.
If $p$ denotes the least prime factor of $n$,
then the expected number of operations required to factor $n$
with the elliptic curve method is
\begin{equation}\label{PreobECCompl}
\exp\left((\sqrt 2+o(1))(\log p)^{1/2}(\log\log p)^{1/2}\right)
(\log n)^{C_2},
\end{equation}
where $C_2$ is a positive constant.
The number field sieve was suggested in~[5].
The heuristic complexity for factoring $n$ via this method is
$$
\exp\left((C_1+o(1))(\log n)^{1/3}(\log\log n)^{2/3}\right).
$$
The main idea underlying the approach suggested in this note is to compute
Fourier coefficients of modular forms that arise from elliptic curves
(see~[6, 7]).

Here is an outline of the method.
It is known that
certain modular forms have integer Fourier coefficients $c_n$,
and these coefficients satisfy specific equations
(see below).
Moreover, the famous Shimura--Taniyama conjecture
(proved by Wiles for the case of
semistable elliptic curves) states that under certain conditions
the numbers $a_p$ --- defined in terms of the orders over ${\mathbb Z}_p$
of a fixed elliptic curve $E$ over ${\mathbb Q}$ --- coincide with the Fourier coefficients $c_p$
of a modular form of the aforementioned type. The numbers $a_p$ are defined by the equation:
$\#E({\mathbb Z}_p)=p+1-a_p$. One may use these properties
to factor $n=pq$. Let $p$ and $q$ be
odd primes of about the same size (such $n$'s
are called RSA-numbers). Choose arbitrary elliptic curve
$E$ over ${\mathbb Q}$. If for the corresponding modular form
we could find the coefficients $c_n$ and $c_{n^2}$,
factor them, and thereby find, say, $c_p$ and $c_{p^2}$,
then from the aforementioned equations for the coefficients
we would find $p={(c_p)}^2-c_{p^2}$.
The problem is to find $c_n$ and $c_{n^2}$ when
the factorization $n=pq$ is unknown. (If the factorization was known,
we would find $c_n$ and $c_{n^2}$ by computing $c_p=a_p$ and $c_q=a_q$
via the Schoof algorithm for the curve $E$.) Thus,
we wish to reduce the problem of factoring $n$
to the problem of factoring $c_n$ and $c_{n^2}$
that depend on a randomly chosen curve $E$.

Assuming that the modular form on the interval $[-1,1]$ has sufficiently
good ap\-pro\-xi\-ma\-ti\-on by partial sums of the Legendre polynomial expansion,
one could try to find $c_n$ and $c_{n^2}$,
factorization of $n$ being unknown, via the ``ap\-pro\-xi\-ma\-ti\-on''
method as it is enunciated at the beginning of the section. One represents
Fourier coefficients of the modular form as integrals, uses the polynomial approximation
and the Schoof algorithm for computing
small sets of coefficients $c_{p'}$ and $c_{p''}$ for small sets of primes
$\{p'\}$ and $\{p''\}$, near $n$ and $n^2$, respectively.
The coefficients in the ``knots'', the sets $\{p'\}$ and $\{p''\}$,
``interpolate'' the coefficients in the ``points'' $n$ and $n^2$.

\medskip

{\bf 2. Main theorem.}
Suppose $\tau=\rho+i\sigma$, $\sigma=\frac1n$ is fixed,
$\rho\in[-1,1]$. Write the function $\phi(\tau)=f(\tau)e^{-2\pi in\tau}$
as the sum of its real and imaginary part\textup{:}
\[
\phi(\tau)=f(\tau)e^{-2\pi in\tau}=u_{\sigma}(\rho)+iv_{\sigma}(\rho).
\]
Suppose
\[
n<p_1'=n+\Delta_1<p_2'=n+\Delta_2<{\ldots}<p_K'=n+\Delta_K.
\]
For an integer $\Delta$, denote by $c_{\Delta}(f)$
the $\Delta$-th Fourier coefficient of $\phi(\tau)$.
It coincides with $(n+\Delta)$-th coefficient of $f(\tau)$\textup{:}
\[
2c_{\Delta}(f) = e^{2\pi\Delta/n}\inte_{-1}^1\(u_{\sigma}(\rho)+iv_{\sigma}(\rho)\)
e^{-i(2\pi\Delta)\rho}\,d\rho.
\]
Take the finite Legendre polynomial expansion
for the real and the imaginary part of $\phi(\tau)$\textup{:}
\be\label{PreobUKVK}
\begin{split}
u_{\sigma}(\rho) &\sim \al_1P_1(\rho)+\ldots+\al_KP_K(\rho)=U_K(\rho),\\
v_{\sigma}(\rho) &\sim \bt_1P_1(\rho)+\ldots+\bt_KP_K(\rho)=V_K(\rho).
\end{split}
\ee
Denote by $C_{\Delta}(\baa,\bb)$ the approximation
to the Fourier coefficient $c_{\Delta}(f)$
obtained via the expansion~\eqref{PreobUKVK}\textup{:}
\[
2C_{\Delta}(\baa,\bb) = e^{2\pi\Delta/n}\inte_{-1}^1\(U_K(\rho)+iV_K(\rho)\)
e^{-i(2\pi\Delta)\rho}\,d\rho.
\]
Making successive substitutions $\Delta=\Delta_1$, $\Delta=\Delta_2$, ${\ldots}$,
$\Delta=\Delta_K$ and integrations, we have the matrix equation
(here $e^{2\pi\Delta/n}$ denotes the diagonal matrix
with the entries $e^{2\pi\Delta_1/n}$, ${\ldots}$, $e^{2\pi\Delta_K/n}$)
\begin{equation*}
\begin{pmatrix}
2C_{\Delta_1}(\baa,\bb) \\
\vdots \\
2C_{\Delta_K}(\baa,\bb)
\end{pmatrix}
=e^{2\pi\Delta/n}
\begin{pmatrix}
a_{1,1}(\Delta_1), & \dots, & a_{1,K}(\Delta_1) \\
\hdotsfor{3} \\
a_{K,1}(\Delta_K), & \dots, & a_{K,K}(\Delta_K)
\end{pmatrix}
\begin{pmatrix}
\al_1+i\bt_1 \\
\vdots \\
\al_K+i\bt_K
\end{pmatrix},
\end{equation*}
or
\begin{equation*}
\begin{pmatrix}
2C_{\Delta_1}(\baa,\bb) \\
\vdots \\
2C_{\Delta_K}(\baa,\bb)
\end{pmatrix}
=e^{2\pi\Delta/n}\BA_{K\times K}
\begin{pmatrix}
\al_1+i\bt_1 \\
\vdots \\
\al_K+i\bt_K
\end{pmatrix}.
\end{equation*}
The entries of the matrix $\BA_{K\times K}$ are Gegenbauer's integrals,
a generalization of Poisson's integral. These integrals can be cast in terms of Bessel functions
of half-integer order, and we have the general formula:
\[
a_{k,m}(\Delta_k) =
(-i)^m\frac1{\sqrt{\Delta_k}}J_{m+1/2}(2\pi\Delta_k).
\]
\begin{theorem}\label{PreobRFC}
The determinant $\det\BA_{K\times K}$ is well-defined when $\Delta_k\ne0$
and nonzero if $\Delta_i\ne\Delta_j$ $(i\ne j)$.
If we define Hankel symbols by
\be\label{PreobHankel}
(\nu,0)=1,\ \ (\nu,l)=\frac{\(4\nu^2-1^2\)\(4\nu^2-3^2\)\ldots\(4\nu^2-(2l-1)^2\)}{2^{2l}l!}\ \
\text{for}\ \ l=1,2,3,\ldots,
\ee
then we can write the formula\textup{:}
\be\label{PreobdetBA}
\begin{split}
\det\BA_{K\times K} = \frac1{\pi^K}&\left(\prod_{\mu=1}^{\flr{(K+1)/2}}i(-1)^{\mu-1}
\frac{(2\mu-1/2,2\mu-2)}{(4\pi)^{2\mu-2}}\right)\times{}\\
{} \times &\left(\prod_{\mu=1}^{\flr{K/2}}(-1)^{\mu-1}
\frac{(2\mu+1/2,2\mu-1)}{(4\pi)^{2\mu-1}}\right)
\left(\prod_{k=1}^K\frac1{\Delta_k}\right)
\left(\prod_{1\le i<j\le K}\left(\frac1{\Delta_j}-\frac1{\Delta_i}\right)\right).
\end{split}
\ee
Define $D_{K,\Delta_0}(E)$ when $\Delta_0\ne0$ to be the supremum of the set
$\{D_{K,\Delta_0}(E)\}$ of numbers satisfying the condition\textup{:}
if $\max(|\e_1'+i\e_1''|,\ldots,|\e_K'+i\e_K''|)<D_{K,\Delta_0}(E)$ then we have
\begin{equation*}
\BA^{-1}_{K\times K}e^{-2\pi\Delta/n}
\begin{pmatrix}
2C_{\Delta_1}+\e_1'+i\e_1'' \\
\vdots \\
2C_{\Delta_K}+\e_K'+i\e_K''
\end{pmatrix}
=
\begin{pmatrix}
\al_1+i\bt_1+\de_1'+i\de_1'' \\
\vdots \\
\al_K+i\bt_K+\de_K'+i\de_K''
\end{pmatrix}
\end{equation*}
and
\[
|2C_{\Delta_0}(\baa+\bde',\bb+\bde'') - 2c_{\Delta_0}(f)|<E.
\]
It follows that if
\be\label{PreobRKf}
R_K(f)=\sum_{k=K+1}^{\infty}\frac2{2k+1}\(\al_k^2+\bt_k^2\)<
\frac1{2\max_{k}e^{4\pi\Delta_k/n}}\(D_{K,\Delta_0}(E)\)^2,
\ee
then
\begin{equation}\label{Preob2cDeltakalphakbetak}
\BA^{-1}_{K\times K}e^{-2\pi\Delta/n}
\begin{pmatrix}
2c_{\Delta_1}(f) \\
\vdots \\
2c_{\Delta_K}(f)
\end{pmatrix}
=
\begin{pmatrix}
\al_1+i\bt_1+\de_1'+i\de_1'' \\
\vdots \\
\al_K+i\bt_K+\de_K'+i\de_K''
\end{pmatrix}
\end{equation}
and
\be\label{Preob2CDelta02cDelta0}
|2C_{\Delta_0}(\baa+\bde',\bb+\bde'') - 2c_{\Delta_0}(f)|<E.
\ee
So if we know the Fourier coefficients
$2c_{\Delta_1}(f)$, ${\ldots}$, $2c_{\Delta_K}(f)$,
then we may recover the Fourier coefficient $2c_{\Delta_0}(f)$ $(\Delta_0\ne0)$
to the precision $E$.
\end{theorem}

\textbf{Remark.}
In this theorem, instead of using Legendre polynomials,
one could use the simple powers $x$, $x^2$, ${\ldots}$, $x^K$.
But in practice, we cannot know that~\eqref{PreobRKf} holds,
and with Legendre polynomials we may infer this from
stabilization of the coefficients
\[
\begin{pmatrix}
\al_1+i\bt_1+\de_1'+i\de_1'' \\
\vdots \\
\al_K+i\bt_K+\de_K'+i\de_K''
\end{pmatrix}
\]
in~\eqref{Preob2cDeltakalphakbetak} as $K$ grows.
For simple powers such stabilization of coefficients
may not occur.

\dokvo We now prove~\eqref{PreobdetBA}. Using the notation~\eqref{PreobHankel}
and known Hankel expansions for Bessel functions of half-integer order,
we deduce that $a_{k,m}(\Delta_k)$ are polynomials in $\frac1{\Delta_k}$
of degree $m$ with the leading coefficients
\[
b_{k,m}=\begin{cases}
\frac{(-i)^m}{\pi}\cos\left(-\frac{(m+1/2)\pi}2-\frac{\pi}4\right)(-1)^{\mu-1}
\frac{(m+1/2,m-1)}{(4\pi)^{m-1}}, &\text{if $m=2\mu-1$};\\
\frac{(-i)^m}{\pi}(-1)\sin\left(-\frac{(m+1/2)\pi}2-\frac{\pi}4\right)(-1)^{\mu-1}
\frac{(m+1/2,m-1)}{(4\pi)^{m-1}}, &\text{if $m=2\mu$}.
\end{cases}
\]
By elementary operations, the matrix $\BA_{K\times K}$ may be transformed
to the Vandermonde matrix.
This completes the proof of~\eqref{PreobdetBA}.

We now prove that the condition~\eqref{PreobRKf}
implies~\eqref{Preob2cDeltakalphakbetak} and~\eqref{Preob2CDelta02cDelta0}.
Define
\[
\phi_{\sigma}(\rho)=u_{\sigma}(\rho)+iv_{\sigma}(\rho),\quad
\Phi_K(\rho)=U_K(\rho)+iV_K(\rho).
\]
We have
\[
\begin{split}
&\inte_{-1}^1|\phi_{\sigma}(\rho)-\Phi_K(\rho)|^2\,d\rho=
\inte_{-1}^1(\phi_{\sigma}(\rho)-\Phi_K(\rho))
\(\overline{\phi_{\sigma}}(\rho)-\overline{\Phi_K}(\rho)\)\,d\rho={} \\
{} = &\inte_{-1}^1|\phi_{\sigma}(\rho)|^2\,d\rho-
\inte_{-1}^1\Phi_K(\rho)\overline{\phi_{\sigma}}(\rho)\,d\rho-
\inte_{-1}^1\overline{\Phi_K}(\rho)\phi_{\sigma}(\rho)\,d\rho+
\inte_{-1}^1|\Phi_K(\rho)|^2\,d\rho={} \\
{} = &\inte_{-1}^1|\phi_{\sigma}(\rho)|^2\,d\rho-
\sum_{k=1}^{K}\frac2{2k+1}\(\al_k^2+\bt_k^2\)=
\sum_{k=K+1}^{\infty}\frac2{2k+1}\(\al_k^2+\bt_k^2\)=R_K(f).
\end{split}
\]
Let us estimate $|2C_{\Delta_k} - 2c_{\Delta_k}(f)|$.
Using the Cauchy inequality and the condition~\eqref{PreobRKf}, for $1\le k\le K$ we obtain the estimates
\[
\begin{split}
|2C_{\Delta_k} - 2c_{\Delta_k}(f)| &\le e^{2\pi\Delta_k/n}
\inte_{-1}^1|\phi_{\sigma}(\rho)-\Phi_K(\rho)|\,d\rho\le{} \\
{} &\le e^{2\pi\Delta_k/n}\sqrt2
\left(\inte_{-1}^1|\phi_{\sigma}(\rho)-\Phi_K(\rho)|^2\,d\rho\right)^{1/2}<
D_{K,\Delta_0}(E),
\end{split}
\]
from which by the definition of $D_{K,\Delta_0}(E)$
we infer~\eqref{Preob2cDeltakalphakbetak} and~\eqref{Preob2CDelta02cDelta0}.
This completes the proof of the theorem.

\medskip

{\bf 3. Modular forms and modular elliptic curves.}
We follow the book~[8].

An unrestricted modular form of level $N$ and
weight 2 is an analytic function
$f$ on $H$ such that
for all $\left(\!\begin{array}{c} a\quad b \\ c\quad d \end{array}\!\right)\in\Gamma_0(N)$
and $\tau\in H$ we have
\begin{equation*}
f\left(\frac{a\tau+b}{c\tau+d}\right)  = (c\tau+d)^2 f(\tau).
\end{equation*}

Since
$\left(\!\begin{array}{c} 1\quad 1 \\ 0\quad 1 \end{array}\!\right) \in \Gamma_0(N)$
then $f(\tau+1)=f(\tau)$ for each modular form and every $\tau$.
Thus $f$ has a Fourier expansion, which is of the form
\begin{equation}\label{PreobQExp}
f(\tau)\ =\ \sum_{n=0}^{\infty} c_ne^{2\pi in\tau}.
\end{equation}
Here
\begin{equation*}
c_n\ =\ \inte_{-\frac12}^{\frac12}f(\tau)e^{-2\pi in\tau}\,d\rho,
\end{equation*}
where $\tau=\rho+i\sigma$, and $\sigma>0$ is fixed.
For cusp forms, $c_0 = 0$.

The following theorem conveys information about the magnitude of the modulus of a cusp form $f(\tau)$
and the Fourier coefficients $c_n$.

\textbf{Theorem ({\rm Deuring}).}
\emph{Let $f\in{\cal S}_2(N)$ have $q$-expansion~{\rm(\ref{PreobQExp})}
at the cusp $\infty$. Then the function $\phi(\tau)=|f(\tau)|\sigma$
is bounded on $H$ and invariant under $\Gamma_0(N)$.
Furthermore, we have $|c_n|\leqslant Cn$.}

Suppose $f\in{\cal S}_2(N)$ is an eigenform, normalized so that
the $q$ expansion~{\rm(\ref{PreobQExp})} has $c_1=1$.
Then for $r\geqslant1$ the Fourier coefficients $c_n$
of~$f$ satisfy
\begin{alignat}{2}
&c_{p^r}c_p=c_{p^{r+1}}+pc_{p^{r-1}},&\quad &\text{for $p$ prime},\ p\nmid N;
\label{PreobEFCoeff1} \\
&c_{p^r}=(c_p)^r,&\quad &\text{for $p$ prime},\ p\mid N; \notag \\
&c_mc_n=c_{mn},&\quad &\gcd(m,n)=1. \notag
\end{alignat}

The following statement is the Shimura--Taniyama conjecture,
proved by Wiles for the case of
semistable elliptic curves
(i.e. for squarefree $N$).

\textbf{Proposition.}
\emph{If $E$ is any elliptic curve defined over ${\mathbb Q}$, if $N$~is its conductor, then
there is a normalized new cusp eigenform $f$ of level $N$ and weight $2$, whose Fourier coefficients
$c_n$ are integers and such that for every prime~$p$ not dividing $N$ $c_p=a_p$ {\rm(}where $a_p$
is defined by $\#E({\mathbb F}_p)=p+1-a_p$,
$\#E({\mathbb F}_p)$ being the order of the group of the elliptic curve $E$ over
${\mathbb F}_p${\rm)}.}

Now we give an example of a modular form associated with an elliptic curve.

{\bf Example ({\rm Hecke}).} Define
\[
\Delta(\tau)=(2\pi)^{12}q\prod_{n=1}^{\infty}\left(1-q^n\right)^{24}
\text{ and }
\eta(\tau)=\frac1{\sqrt{2\pi}}\Delta(\tau)^{\frac1{24}}.
\]
Then it follows that $f(\tau)=\eta(11\tau)^2\eta(\tau)^2$
is a new cusp eigenform of weight $2$ and level $11$.
It is associated with the elliptic curve $E$
\[
y^2+y=x^3-x^2-10x-20,
\]
which for $p\ne2,3$ can be given the form
\[
y^2=x^3-\frac{31}3x-\frac{2431}{108}.
\]

{\bf 4. Elliptic curves over finite fields.}
The following theorem gives a range for the order of an elliptic curve group
defined over a finite field.

\textbf{Theorem ({\rm Hasse}).}
\emph{Let $E$ be an elliptic curve over ${\mathbb Q}$ with integer coefficients
and the discriminant $\Delta$.
For each prime~$p\nmid\Delta$, let $E({\mathbb F}_p)$ be the
reduction of $E$ modulo $p$.
Then $|p+1-\#E({\mathbb F}_p)|<2\sqrt p$.}

We need two celebrated algorithms related to
elliptic curves over finite fields:
the Schoof algorithm~[9] for elliptic curve point-counting
and Lenstra's elliptic curve method
for factorization, mentioned in Section 1.
For prime $q>3$ and elliptic curve $E({\mathbb F}_q)$,
the Schoof algorithm computes $\#E({\mathbb F}_q)$ in complexity $O\left({\log}^8q\right)$.
The Lenstra method factors a composite $n$.
If $p$ denotes the least prime factor of $n$,
then this method has the complexity~(\ref{PreobECCompl}).

\medskip

{\bf 5. Smooth numbers.} A positive integer is said to be $y$-smooth
if it does not have any prime factor exceeding $y$.
Let $\psi(x,y)$ denote the number of $y$-smooth integers
among positive integers $n\leqslant x$.

\textbf{Theorem ({\rm Canfield--Erd\H os--Pomerance [10]}).}
\emph{The estimate $\psi(x,x^{1/u})=xu^{-u+o(u)}$
holds uniformly
as $u\to\infty$, $u<(1-\varepsilon_1)\log x/\log\log x$, $\varepsilon_1>0$.}

\medskip

{\bf 6. The central result and the algorithm.}
\begin{theorem}
Assuming that the method discussed in Section $2$ for computing the Fourier coefficients
of a modular form associated with an elliptic curve is correct
and runs in polynomial time, the heuristic complexity estimate
for factoring an integer $n$
with the least prime factor $p$ {\rm(}$p^2\nmid n${\rm)} is
\[
\exp\left(c_1(k)(\log p)^{1/k}(\log\log p)^{1-1/k}\right)(\log n)^{c_2(k)}
\]
arithmetic operations for any fixed positive integer $k${\rm;}
here $c_1(k)$, $c_2(k)$ are positive constants which may depend
on the choice of $k$.
\end{theorem}

We now describe a recursive algorithm having the claimed complexity estimate.
First, we show that the following method gives the complexity estimate
\[
\exp\left(c_1(3)(\log p)^{1/3}(\log\log p)^{2/3}\right)(\log n)^{c_2(3)}
\]
(this method corresponds to the case $k=3$).

Step 0. Choose a smoothness parameter
$B\asymp p^{1/\left((\log p)^{1/3}(\log\log p)^{-1/3}\right)}$, where
$p$ is the least prime factor of $n$
to be factored.
(In fact, since we do not know $p$ to begin with,
we are instructed to start with a low $B_0$ value
and then to run the algorithm with $B=B_0$, $B=2B_0$, $B=4B_0$, $\ldots$,
choosing for each $B$ about
$$
\exp\left((\sqrt2+o(1))(\log B)^{1/2}(\log\log B)^{1/2}\right)
$$
random elliptic curves in Step~1.)

Step 1. Choose random elliptic curve $E$ over ${\mathbb Q}$
with integer coefficients.

Step 2. Compute the Fourier coefficients $c_n$ and $c_{n^2}$
of the associated modular form with the method discussed in Section 2.

Step 3. Try to factor these coefficients
with the elliptic curve method, hoping for
a successful event:
the largest prime factor of the coefficient $c_p\mid c_n$
and the largest prime factor of the coefficient $c_{p^2}\mid c_{n^2}$ both do not exceed $B$.
For each $c_1\mid c_n$ and $c_2\mid c_{n^2}$
compute $d=\gcd((c_1)^2-c_2,n)$ and, if $1<d<n$,
return the nontrivial divisor $d$ of $n$.

Step 4. Failure: goto Step~$1$ (or give up on the factorization attempt).

Assume that the Fourier coefficients $c_p$ and $c_{p^2}$ of the modular form
associated with the chosen elliptic curve
(where $c_p=O\left(\sqrt p\right)$ and $c_{p^2}=O(p)$ by the Hasse theorem
and the relation~(\ref{PreobEFCoeff1}))
are $B$-smooth with the same probability
as random integers chosen from the respective intervals.
Then we get $B$-smooth coefficients $c_p$ and $c_{p^2}$
expecting about $u^{\frac32u}$, where $u=(\log p)^{1/3}(\log\log p)^{-1/3}$,
elliptic curves, with the attempt to find these coefficients
in Step~$3$ requiring
$$
\exp\left((\sqrt2+o(1))(\log B)^{1/2}(\log\log B)^{1/2}\right)(\log n)^{C_2'}
$$
arithmetic operations.
Multiplying the expressions, we obtain the estimate for the complexity of the algorithm in the form
\[
\exp\left(c_1(3)(\log p)^{1/3}(\log\log p)^{2/3}\right)(\log n)^{c_2(3)}.
\]

If we now use this algorithm instead of the elliptic curve method in Step~$3$
and choose the optimal smoothness parameter,
then we get an algorithm having complexity
\[
\exp\left(c_1(4)(\log p)^{1/4}(\log\log p)^{3/4}\right)(\log n)^{c_2(4)}.
\]

Furthermore, if we use $k+1$ levels of recursion, choose
$
B\asymp p^{1/\left((\log p)^{1/(k+1)}(\log\log p)^{-1/(k+1)}\right)}
$
and apply the algorithm with $k$ levels of recursion,
having complexity
\[
\exp\left(c_1(k)(\log p)^{1/k}(\log\log p)^{1-1/k}\right)(\log n)^{c_2(k)},
\]
then we arrive at the following complexity estimate:
\begin{gather*}
\left((\log p)^{1/(k+1)}(\log\log p)^{-1/(k+1)}\right)^{\frac32
(\log p)^{1/(k+1)}(\log\log p)^{-1/(k+1)}}\,{\times}\\
{\times}\,\exp\left(c(k)(\log B)^{1/k}(\log\log B)^{1-1/k}\right)(\log n)^{c_2'(k)}\,{=}\\
{=}\,\exp\left(c_1'(k)(\log p)^{1/(k+1)}(\log\log p)^{1-1/(k+1)}\right)\,{\times}\\
{\times}\,\exp\left(c_1''(k)(\log p)^{(1/k)(1-1/(k+1))}(\log\log p)^{1/(k(k+1))+1-1/k}\right)(\log n)^{c_2'(k)}\,{=}\\
{=}\,\exp\left(c_1((k+1))(\log p)^{1/(k+1)}(\log\log p)^{1-1/(k+1)}\right)(\log n)^{c_2((k+1))}.
\end{gather*}

\makeatletter
\def\@biblabel#1{#1.}
\makeatother


\end{document}